\documentclass[10pt,twoside,reqno,a4paper]{amsart}
\usepackage[utf8]{inputenc}
\usepackage{courier}
\usepackage[T1]{fontenc}

\usepackage[mathbf]{euler}
\usepackage{topformat}
\usepackage[driver=pdftex,margin=3cm,heightrounded=true,centering]{geometry}
\usepackage{hyperref}
\usepackage{enumerate}
\usepackage{amssymb}
\usepackage{amsopn}
\usepackage{amsmath}
\usepackage{mathtools}
\usepackage{tikz}
\usepackage{tikz-cd}
\usepackage{ifthen}
\usepackage{graphics}
\usepackage{amsthm}
\usepackage{topthm}
\author{Thorben Kastenholz}
\thanks{This research is funded by the Deutsche Forschungsgemeinschaft (DFG,
German Research Foundation) – Project number 281869850}
\date{\today}
\title{Simplicial volume of open books in dimension $4$}
\address{Karlsruher Institut f\"ur Technologie, Englerstraße 2, 76131
  Karlsruhe, Germany}
\email{thorben.kastenholz@kit.edu}
\begin{document}
\newcommand{\introduce}[1]
  {\textbf{#1}}
\newcommand{\tk}[1]{\todo[size=\tiny,color=green!40]{TK: #1}}
\newcommand\blfootnote[1]{%
  \begingroup
  \renewcommand\thefootnote{}\footnote{#1}%
  \addtocounter{footnote}{-1}%
  \endgroup
}

\newcommand{\apply}[2]
  {{#1}\!\left({#2}\right)}
\newcommand{\at}[2]
  {\left.{#1}\right\rvert_{#2}}
\newcommand{\Identity}%
  {\mathrm{Id}}
\newcommand{\NaturalNumbers}%
  {\mathbf{N}}
\newcommand{\Integers}%
  {\mathbf{Z}}
\newcommand{\Rationals}%
  {\mathbf{Q}}
\newcommand{\Reals}%
  {\mathbf{R}}
  \newcommand{\ComplexNumbers}%
  {\mathbf{C}}
\newcommand{\AbstractProjection}[1] 
  {p_{#1}}
\newcommand{\RealPart}[1]
  {\apply{\operatorname{Re}}{#1}}
\newcommand{\ImaginaryPart}[1]
  {\apply{\operatorname{Im}}{#1}}
\newcommand{\Floor}[1]
  {\left \lfloor #1 \right \rfloor}
\newcommand{\Norm}[1]
  {\left|\left|#1\right|\right|}

\newcommand{\Surface}[1]
  {\Sigma_{#1}}
\newcommand{\SurfaceGroup}[1]
  {S_{#1}}
\newcommand{\Manifold}%
  {M}
\newcommand{\Page}%
  {P}
\newcommand{\Binding}%
  {B}
\newcommand{\OpenBook}[2]
  {\apply{o}{#1,#2}}
\newcommand{\Tangenbundle}[1]
  {T#1}
\newcommand{\NormalBundle}[2]
  {N_{#2}#1}
\newcommand{\FiberTransferHomology}[1]
  {#1^{!}}
\newcommand{\Submanifold}%
  {I}
\newcommand{\TubularNeighborhood}[1]
  {U_{#1}}
\newcommand{\ManifoldAlternative}%
  {N}
 \newcommand{\ManifoldAuxiliary}%
  {K}
\newcommand{\NullBordism}%
  {W}
\newcommand{\Bordism}
  {P}
 \newcommand{\ManifoldFiber}%
 {M}
 \newcommand{\ManifoldTotal}%
 {E}
  \newcommand{\ManifoldBase}%
  {B}
\newcommand{\SmoothMap}%
  {\phi}
  \newcommand{\MorseFunction}%
  {f}
\newcommand{\Diffeomorphism}%
  {\Phi}
\newcommand{\Dimension}%
  {d}
  \newcommand{\HalfDimension}%
  {n}
\newcommand{\FundamentalClass}[1]
  {\left[#1\right]}
\newcommand{\Interval}%
  {I}
\newcommand{\Ball}[1]
  {D^{#1}}
\newcommand{\Sphere}[1] 
  {S^{#1}}
\newcommand{\Torus}[1]
  {T^{#1}}
\newcommand{\SimplicialVolume}[1]
  {\lvert \lvert #1 \rvert \rvert}
\newcommand{\ellone}%
  {\ell_{1}}
\newcommand{\Boundary}[1]
  {\partial #1}
\newcommand{\ComplexOfEmbeddings}[1]
  {\apply{K}{#1}}
\newcommand{\GenusOf}[1]
  {\apply{\Genus}{#1}}
\newcommand{\StableGenusOf}[1]
  {\apply{\overline{\Genus}}{#1}}
\newcommand{\Surgery}%
  {\natural}
\newcommand{\HandleEmbedding}%
  {\Phi}
\newcommand{\Singularity}%
  {X}

\newcommand{\Diff}[1]
  {\mathrm{Diff}\!\left(#1\right)}
\newcommand{\DiffGroup}[1]
  {\mathrm{Diff}^{B\Group}\!\left(#1\right)}
\newcommand{\DiffZero}[1]
  {\mathrm{Diff}_0\!\left(#1\right)}
\newcommand{\DiffOne}[1]
  {\widetilde{\mathrm{Diff}}_0\!\left(#1\right)}
\newcommand{\HomeoGroup}[1]
  {\apply{\mathrm{Homeo}}{#1}}
\newcommand{\HomeoCompactGroup}[1]
  {\apply{\mathrm{Homeo}_{c}}{#1}}
\newcommand{\HomeoLowerGroup}[1]
  {\apply{\mathrm{Homeo}^{\geq}}{#1}}

\newcommand{\FiberingSpace}%
  {E}
\newcommand{\FiberingProjektion}[1] 
  {\pi_{#1}}
\newcommand{\Fiber}%
  {F}
\newcommand{\FiberDimension}%
  {f}
\newcommand{\Base}%
  {B}
\newcommand{\ClutchingFunction}[1] 
  {\varphi_{#1}}

\newcommand{\Group}%
  {\Gamma}
\newcommand{\Subgroup}
  {H}
\newcommand{\IntegralGroupRing}[1]
  {\Integers#1}
\newcommand{\AmenableGroup}
  {A}
\newcommand{\GroupElement}%
  {g}
\newcommand{\Genus}%
  {g}
\newcommand{\QuadraticModule}%
  {\mathbf{M}}
\newcommand{\WittIndex}[1]
  {\apply{\Genus}{#1}}
\newcommand{\StableWittIndex}[1]
  {\apply{\overline{\Genus}}{#1}}
\newcommand{\ComplexOfHyperbolicInclusions}[1]
  {\apply{K^{a}}{#1}}
\newcommand{\ChainContraction}[1]
  {H_{#1}}
\newcommand{\FreeGroup}[1]
  {F_{#1}}

\newcommand{\HomologyClass}%
  {\alpha}
\newcommand{\HomologyOfSpaceObject}[3]
  {\apply{H_{#1}}{#2 ; #3}}
\newcommand{\CohomologyOfSpaceObject}[3]
  {\apply{H^{#1}}{#2 ; #3}}
\newcommand{\BoundedCohomologyOfSpaceObject}[3]
  {\apply{H^{#1}_{\text{b}}}{#2 ; #3}}
\newcommand{\BoundedCohomologyOfSimplicialObject}[3]
  {\apply{H^{#1}_{\text{b, s}}}{#2 ; #3}}
\newcommand{\HomologyOfSpaceMorphism}[1]
  {{#1}_{\ast}}
\newcommand{\HomologyOfGroupObject}[3]
  {\apply{H_{#1}}{#2; #3}}
\newcommand{\HomologyOfGroupMorphism}[2]
  {{#1}_{\ast}}
\newcommand{\HomologyOfSpacePairObject}[3]
  {\apply{H_{#1}}{{#2},{#3}}}
\newcommand{\Multiple}%
  {\lambda}

\newcommand{\TopologicalSpace}%
  {X}
\newcommand{\MappingTorus}[1]
  {T_{#1}}
\newcommand{\Inner}[1]
  {\mathring{#1}}
\newcommand{\Point}%
  {\ast}
\newcommand{\Loop}%
  {\gamma}
\newcommand{\ContinuousMap}%
  {f}
  \newcommand{\ContinuousMapALT}%
  {g}
 \newcommand{\maps}%
  {\ensuremath{\text{maps}}}
\newcommand{\HomotopyGroupOfObject}[3]
  {\apply{\pi_{#1}}{{#2},{#3}}}
\newcommand{\HomotopyGroupOfPairObject}[4]
  {\apply{\pi_{#1}}{{#2},{#3},{#4}}}
\newcommand{\HomotopyGroupMorphism}[1] 
  {{#1}_{\ast}}
\newcommand{\EMSpace}[2]
  {\apply{K}{{#1},{#2}}}
\newcommand{\ClassifyingSpace}[1] 
  {B#1}
\newcommand{\UniversalCovering}[1] 
  {\widetilde{#1}}
\newcommand{\UniversalCoveringMap}[1] 
  {\widetilde{#1}}
\newcommand{\FundamentalCycle}[1]
  {\sigma_{#1}}

\newcommand{\SimplicialComplex}
  {X}
\newcommand{\AuxSimplicialComplex}
  {K}
\newcommand{\Subcomplex}
  {Y}
\newcommand{\AuxSubcomplex}
  {L}
\newcommand{\Simplex}[1]
  {\sigma_{#1}}
\newcommand{\Link}[2]
  {\apply{\text{Lk}_{#1}}{#2}}
\newcommand{\Star}[2]
  {\apply{\text{St}_{#1}}{#2}}
\newcommand{\BoundaryIndexSimplex}[2]
  {\apply{\partial_{#1}}{#2}}
\newcommand{\BoundarySimplex}%
  {\partial}
\newcommand{\StandardSimplex}[1]
  {\Delta_{#1}}
\newcommand{\vertex}
  {v}
\newcommand{\GeometricRealization}[1]
  {\left\lvert #1 \right\rvert}
\newcommand{\BoundHomotopy}
  {N}
\newcommand{\Horn}[2]
  {\Lambda^{#1}_{#2}}

\begin{abstract}
  In this short note we adapt a proof by Bucher and Neofytidis to prove that
  the simplicial volume of $4$-manifolds admitting an
  open book decomposition vanishes. In particular this shows that Quinns
  signature invariant, which detects the existence of an open book
  decomposition in dimensions above $5$, is insufficient to
  characterize open books in dimension $4$, even if one allows arbitrary
  stabilizations via connected sums.
\end{abstract}

\maketitle
\section{Introduction}
Open book decompositions are very useful decompositions of manifolds. They are
most famously used in dimension $3$, but also useful in different dimensions.
An $n$-dimensional open book consists of an $n-2$-dimensional manifold
$\Binding$ called the \introduce{binding}, an $n-1$-dimensional manifold
$\Page$ with boundary equal to $\Binding$ called the \introduce{page} and
a
diffeomorphism $\Diffeomorphism\colon
\Page \to \Page$ that restricts to the identity on the boundary called
the \introduce{monodromy}.
We say that an $n$-manifold $\Manifold$ admits an open book decomposition if
there exists an $n$-dimensional open book such that $\Manifold$ is
isomorphic to
\[
  \OpenBook{\Page}{\Diffeomorphism}
  \coloneqq
  \MappingTorus{\Diffeomorphism}
  \cup_{\Sphere{1}\times \Binding}
  \Ball{2}\times \Binding
\]%
where $\MappingTorus{\Diffeomorphism}$ denotes the mapping torus of
$\Diffeomorphism$.

The classification of surfaces is sufficient to conclude that the only surfaces
admitting an open book decomposition are the sphere and the torus. It is a
classical result that every $3$-manifold admits an open book decomposition.

In all dimension bigger or equal than $5$, it was established by Quinn in
\cite{Quinn} that the existence of an
open book decompositions for a compact manifold is equivalent to the vanishing
of an algebraic invariant, the so called asymmetric signature (See
\cite{Ranicki} for a definition).
This invariant vanishes in odd dimensions, in particular every odd dimensional
manifold admits an open book decomposition. The proof in \cite{Quinn} relies on
the s-cobordism theorem and hence breaks down in dimension $4$. We show the
following:
\begin{theorem}
\label{thm:NonExistence}
  There exists a $4$-manifold $X$, whose asymmetric signature vanishes, but
  which does not admit an open book decomposition. Moreover for any
  $4$-manifold $Y$, the connected sum $X\# Y$ does not admit an open book
  decomposition as well.
\end{theorem}
The manifold $X$ in \autoref{thm:NonExistence} will be a product of two
surfaces of genus greater or equal than $2$. Note that it follows from
Proposition~30.6 in \cite{Ranicki} that the asymmetric signature is also the
obstruction to the existence of a twisted double structure on a manifold. Since
surfaces are all doubles, it follows that the product of two surfaces has
vanishing asymmetric signature. We will establish the following which also
shows that there are many manifolds in dimension $4$ that do not admit open
book decompositions:
\begin{theorem}
\label{prp:Contradiction}
  The simplicial volume of all $4$-manifold admitting an open book
  decomposition vanishes.
\end{theorem}
Since the product of two surfaces of genus greater or equal than two has
non-zero simplicial volume and the simplicial volume is additive with respect
to connected sums (See \cite{Frigerio} for both of these facts), one
immediately obtains \autoref{thm:NonExistence}.

The proof of \autoref{prp:Contradiction} is based on the proof in
\cite{BucherNeofytidis}, where it is shown that the mapping torus of any closed
$3$-manifold has vanishing simplicial volume. The main reason for their result
as well as for \autoref{prp:Contradiction} lies in the rigidness of
$3$-manifolds. More concretely both proofs use that hyperbolic
$3$-manifolds admit only finitely many diffeomorphisms up to isotopy, while
Seifert fibered $3$-manifolds can not contribute to the simplicial volume of a
mapping torus. The main difference in the proofs is that we do not use the
Geometrization Theorem and that we have to use actual fundamental cycles in
order to glue a fundamental cycle of the mapping torus of the page to a
fundamental cycle of a neighborhood of the binding. However our proof still
follows theirs very closely.

The proof first reduces to the case of an irreducible page and then uses the JSJ
decomposition of an irreducible page to construct arbitrarily small fundamental
cycles piece by piece. It is crucial for this second step that the monodromy
fixes the boundary pointwise.

\paragraph*{Acknowledgments:} The author would like to thank Marc Kegel, Felix
Schm\"aschke and Chun-Sheng Hsueh for many helpful discussion about open book
decompositions in dimension $4$ and for proofreading this note.
\section{Open book decompositions in dimension $4$}
\label{scn:SimpVolume}
\paragraph{Reducible pages and their diffeomorphisms}
Let $\Manifold$ denote a closed $4$ manifold admitting an open book
decomposition with page $\Page$, binding $\Binding= \partial \Page$ and
monodromy $\Diffeomorphism$.

First of all suppose that $\Page$ is a reducible $3$-manifold, then
analogously to the closed case in \cite{BucherNeofytidis} the diffeomorphism
$\Diffeomorphism$ can be written as a composition of first some diffeomorphisms
of the prime summands $\Diffeomorphism_{i}$, then some permutations of
diffeomorphic summands and then some slide homomorphisms. More precisely, by
the Milnor-Kneser theorem $\Page = \#_i \Page_i \# \#_{k} \Sphere{1}\times
\Sphere{2}$, where $\Page_{i}$ is assumed to be irreducible i.e. every embedded
$2$-sphere bounds a ball. In particular every boundary component of $\Page$ is
contained in some $\Page_i$. Then by \cite{McCullough}
any diffeomorphism
$\Diffeomorphism \colon \Page \to \Page$ can be written up to isotopy as a
composition
$\Diffeomorphism_3 \Diffeomorphism_2 \circ \Diffeomorphism_1$, where
\begin{itemize}
\item
  $\Diffeomorphism_1$ is a composition of diffeomorphisms of the prime summands
  $\Page_i$ and the $\Sphere{1}\times \Sphere{2}$ summands,
\item
  $\Diffeomorphism_2$ consists of swaps of diffeomorphic prime summands,
\item
  $\Diffeomorphism_3$ is a composition of slide homomorphisms, where a slide
  homomorphisms is defined as follows: Pick some $i$ or $k$ and let us denote
  by $\overline{\Page}$ the quotient of $\Page$ obtained by collapsing
  $\Page_i$ or the $k$-th $\Sphere{1} \times \Sphere{2}$.
  Let $\alpha$ denote an embedded oriented circle in $\overline{\Page}$ that
  meets the
  quotient point $\Point$ and let $J_t \colon \overline{\Page}\to
  \overline{\Page}$ denote an isotopy that moves $\Point$ along
  $\alpha$. A slide homeomorphism $\phi \colon \Page \to \Page$ is defined by
  $\at{\phi}{\Page \setminus \Page_i} = \at{J_1}{\overline{\Page} \setminus
  \Point}$ and $\at{\phi}{\Page_i} = \Identity$.
\end{itemize}
We will first show that the simplicial volume of $\Manifold$ is dominated by
the sum of the simplicial volumes of all open books corresponding to
the diffeomorphisms occuring in $\Diffeomorphism_1$.

Let $\hat{\Page}$ denote
$\bigvee_i \Page_{i} \vee \bigvee_k \Sphere{1}$ and $\pi \colon \Page
\to \hat{\Page}$, the projection. Analogously to \cite{BucherNeofytidis},
$\Diffeomorphism$ descends to a diffeomorphism $\hat{\Diffeomorphism}$ on
$\hat{\Page}$ defined as follows:
\begin{itemize}
  \item
  $\hat{\Diffeomorphism}$ is defined by descending $\Diffeomorphism_1$,
  $\Diffeomorphism_2$ and $\Diffeomorphism_3$ to $\hat{\Page}$
  \item
  $\hat{\Diffeomorphism}_1$ is given by the same diffeomorphisms of the
  summands,
  \item
  $\hat{\Diffeomorphism}_2$ is given by the corresponding swaps as in
  $\Diffeomorphism_2$
  \item
  for every slide homomorphism that slides $\Page_i$ along the path $\alpha$,
  take a
  ball around the wedge point in $\Page_{i}$, collapse this ball to an interval
  and map this interval along the image $\apply{\pi}{\alpha}$.
  \item
  for every slide homomorphism that slides a $\Sphere{1}\times \Sphere{2}$
  summand
  along the path $\alpha$, map the corresponding $\Sphere{1}$ in $\hat{\Page}$
  to the
  conjugation of itself by $\alpha$.
\end{itemize}
Evidently the diagram
\[
  \begin{tikzcd}
    \Page
      \ar[r,"\Diffeomorphism"]
      \ar[d,"\pi"]
    &
    \Page
      \ar[d,"\pi"]
    \\
    \hat{\Page}
      \ar[r,"\hat{\Diffeomorphism}"]
    &
    \hat{\Page}
  \end{tikzcd}
\]
commutes
Additionally $\pi$ extends to a map
$
  \overline{\pi}
  \colon
  \OpenBook{\Page}{\Diffeomorphism}
  \to
  \OpenBook{\hat{\Page}}{\hat{\Diffeomorphism}}
$%
where $\OpenBook{\hat{\Page}}{\hat{\Diffeomorphism}}$ denotes
$
  \MappingTorus{\hat{\Diffeomorphism}}
  \cup_{\Sphere{1}\times \Binding}
  \Ball{2}\times \Binding
$
. Note that this is possible because every connected component of the boundary
of $\Page$ is contained in a $\Page_{i}$.
Since $\pi$ induces an
isomorphism on fundamental groups and the identity on boundary components,
$\overline{\pi}$ induces an isomorphism on fundamental groups as well.
Therefore the map $\overline{\pi}$ induces an isometry in bounded cohomology
and hence
the simplicial volume of
$\Manifold = \OpenBook{\Page}{\Diffeomorphism}$ agrees with the norm of the
"fundamental class" of
$\OpenBook{\hat{\Page}}{\hat{\Diffeomorphism}}$.

Now define $\Page^{S}$ as
$\hat{\Page}\vee \bigvee_s \Sphere{1}$, where $s$ denotes the number of
slide homomeorphisms in the decomposition of $\Diffeomorphism$. Analogously to
\cite{BucherNeofytidis} we construct a map $\Page^{S} \to \hat{\Page}$ as the
identity on $\hat{\Page} \subset \Page^{S}$ and as the corresponding
$\apply{\pi}{\alpha}$ for every $\Sphere{1}$.
We also define another map
$
  \Diffeomorphism^{S}
  \colon
  \Page^{S}
  \to
  \Page^{S}
$
as follows:
\begin{itemize}
\item
  $\Diffeomorphism^{S}$ is defined by descending $\hat{\Diffeomorphism}_1$,
  $\hat{\Diffeomorphism}_2$ and the descended slide homomorphisms to
  $\Page^{S}$
\item
  $\hat{\Diffeomorphism}_1$ and $\hat{\Diffeomorphism}_2$ are extended via the
  identity from $\hat{\Page} \subset \Page^{S}$ to $\Page^{S}$,
\item
  a slide homomorphism that slides $\Page_i$ along $\alpha$, descends like
  the corresponding part of $\hat{\Diffeomorphism}$ except that it maps the
  collapsed ball to the circle $\Sphere{1} \subset \bigvee_s \Sphere{1} \subset
  \Page^{S}$ corresponding to the slide homomorphism
\item
  a slide homomorphism that slides a circle $\Sphere{1} \subset \bigvee_k
  \Sphere{1} \subset \hat{\Page}$ along a path $\alpha$ is not mapped to the
  conjugation of itself by $\alpha$ but to the conjugation of itself by the
  corresponding loop in $\bigvee_s \Sphere{1} \subset \Page^{S}$.
\end{itemize}
Evidently the diagram
\[
  \begin{tikzcd}
    \Page^{S}
    \ar[r,"\Diffeomorphism^{S}"]
    \ar[d]
    &
    \Page^{S}
    \ar[d]
    \\
    \hat{\Page}
    \ar[r,"\hat{\Diffeomorphism}"]
    &
    \hat{\Page}
  \end{tikzcd}
\]
commutes and furthermore the vertical map sends the "fundamental class" of
$\Page^{S}$ to the "fundamental class" of $\hat{\Page}$.
Therefore the norm of the "fundamental class" of
$\OpenBook{\Page^{S}}{\Diffeomorphism^{S}}$, which again makes sense since the
every boundary components is contained in some $\Page_i$,
is an upper bound for the norm of the "fundamental class" of
$\OpenBook{\hat{\Page}}{\hat{\Diffeomorphism}}$.
By potentially replacing
$\Diffeomorphism^{S}$ by a finite power of itself (which can only increase the
simplicial volume of the corresponding open book since there is a map of
positive degree from the new open book to the old one),
we can assume that
$\Diffeomorphism^{S}$ maps $\Page_{i}$ to $\Page_{i} \vee \bigvee_s
\Sphere{1}$ and restricts to the identity on $\bigvee_s \Sphere{1}$. Hence by
subadditivity of the simplicial volume, we have
\[
  \SimplicialVolume{\OpenBook{\Page^{S}}{\Diffeomorphism^{S}}}
  \leq
  \sum_i
  \SimplicialVolume{
    \OpenBook
      {\Page_{i} \vee \bigvee_s \Sphere{1}}
      {\at
        {\Diffeomorphism^{S}}
        {\Page_{i} \vee \bigvee_s \Sphere{1}}
      }
  }
\]
Now completely analogous to \cite{BucherNeofytidis}, one can
see that if the
simplicial volume of $\OpenBook{\Page_{i}}{\Diffeomorphism_{i}}$
vanishes, then the simplicial volume of
$
  \OpenBook
    {
      \Page_{i} \vee \bigvee_s \Sphere{1}
    }
    {\at
      {\Diffeomorphism^{S}}
      {\Page_{i} \vee \bigvee_s \Sphere{1}}
    }
$
vanishes as well. Indeed, the map $\at{\Diffeomorphism^{S}}{\Page_i}$ collapses
a neighborhood of the wedge point in $\Page_{i}$ to an interval which is then
wrapped around a loop in $\bigvee_s \Sphere{1}$. Hence
$\at{\Diffeomorphism^{S}}{\Page_i}$ is homotopic a diffeomorphism of $\Page_i$
$\tilde{\Diffeomorphism}_i$, hence we obtain a map
$
  \OpenBook{\Page_i}{\tilde{\Diffeomorphism}_i}
  \to
  \OpenBook
    {
      \Page_{i} \vee \bigvee_s \Sphere{1}
    }
    {\at
      {\Diffeomorphism^{S}}
      {\Page_{i} \vee \bigvee_s \Sphere{1}}
    }
$
which induces an isomorphism on $4$-th homology. In particular the simplicial
volume of $\OpenBook{\Page_i}{\tilde{\Diffeomorphism}_i}$ is an upper bound for
the simplicial volume of
$
  \OpenBook
  {
    \Page_{i} \vee \bigvee_s \Sphere{1}
  }
  {\at
    {\Diffeomorphism^{S}}
    {\Page_{i} \vee \bigvee_s \Sphere{1}}
  }
$.
All in all we conclude that the simplicial volume of an open book with
reducible page is bounded from above by the sum of the simplicial volumes of
certain open books with its irreducible components as pages. In particular this
reduces the proof of \autoref{prp:Contradiction} to the case where the page is
irreducible.

\paragraph{Irreducible pages}
We need to better understand the structure of a $3$-manifold with boundary. A
$3$-manifold is called irreducible if every embedded $2$-sphere bounds a ball.
We call it boundary-irreducible if no simple closed curve in the boundary
bounds a $2$-disk in the $3$-manifold. By the Loop Theorem this is equivalent
to the boundary inclusion being componentwise $\pi_1$-injective.
\begin{definition}
  Let $\Submanifold$ denote a possibly disconnected closed surface. A
  \introduce{compression body} is a $3$-manifold $V$, obtained from
  $\Submanifold\times [0,1]$ by attaching $1$-handles to $\Submanifold \times
  \{1\}$. We call $\Submanifold \times \{0\}$ the \introduce{outer boundary of
  $V$} and
  $\partial V \setminus \Submanifold \times \{0\}$ the
  \introduce{inner boundary of $V$}.
\end{definition}
By Theorem~3.7 in \cite{Bonahon} any irreducible $3$-manifold $X$ contains an
up to isotopy unique compression body $V$ such that
\begin{itemize}
\item
  The outer boundary of $V$ consists of the whole boundary of $X$,
\item
  $X \setminus V$ is boundary-irreducible.
\end{itemize}
Similarly to the closed case, any irreducible and boundary-irreducible
$3$-manifold $X$ carries a JSJ-decomposition (See Theorem~3.8 in
\cite{Bonahon}):
There exists an up to isotopy unique collection of tori $S$ and properly
embedded annuli $A$ such that the components of the complement of these tori
and annuli are either hyperbolic or Seifert fibered.

The goal of the remainder of this section is to prove that any open book
with an irreducible page has vanishing simplicial volume. By the previous
subsection, this suffices to prove \autoref{prp:Contradiction}.

Let $\Manifold$ denote an open book with irreducible page $\Page$, binding
$\Binding$ and monodromy $\Diffeomorphism$.
Since $\Page$ is irreducible, let $V$ denote the previously considered
compression body and $S$ and $A$ the collection of tori and annuli of its JSJ
decomposition.
Hence
$\Page$ can be written as the union
\[
  V \cup \bigcup_i \Page_{H,i} \cup \bigcup_j \Page_{S,j}
\]
where $\Page_{H,i}$ denotes the hyperbolic pieces and $\Page_{S,j}$ the Seifert
fibered ones, let us denote the corresponding restrictions of $\Diffeomorphism$
by $\Diffeomorphism_{H,i}$ and $\Diffeomorphism_{S,j}$.

Since the diffeomorphism type of the mapping torus only depends on the isotopy
class of the diffeomorphism, we can assume that $\Diffeomorphism$ maps $V$,
$S$ and $A$ to themselves. Additionally by potentially replacing
$\Diffeomorphism$ by
a power of itself (which can only increase the simplicial volume of the
corresponding open book, since there exists a map of positive degree from the
newly considered open book to the old one), we
can assume that $\Diffeomorphism$ acts as the identity on the set of
connected components of $\Page \setminus (V \cup S \cup A)$ and on the set of
connected components of $S$ and $A$.

We will now proceed by providing representatives of the fundamental class of
$\Manifold$ whose norm will be arbitrarily small. These cycles will be
constructed by defining them on all the
occuring pieces of $\Page$ and showing that they either already agree on their
respective boundaries, or their differing boundaries are supported by a
submanifold with amenable fundamental group and hence can be glued together
efficiently.

Fix representatives $\FundamentalCycle{\Ball{2}}$ and
$\FundamentalCycle{\Binding}$ of the fundamental class of
$\Ball{2}$ and $\Binding$ respectively.
Then
$
  \tau_\Binding
  \coloneqq
  \FundamentalCycle{\Ball{2}}
  \times
  \FundamentalCycle{\Binding}
$
represents the fundamental class of $\Ball{2}\times \Binding \subset
\Manifold$. By choosing $\FundamentalCycle{\Ball{2}}$ sufficiently small, the
norm of $\tau_\Binding$ can be made arbitrarily small.

Let $\FundamentalCycle{V}$ represent the fundamental class of $V$ such that
$\partial \FundamentalCycle{V}$ agrees with $\FundamentalCycle{\Binding}$ on
the corresponding boundary component.
By Lemma~3.2 in \cite{Oertel}, we have that $\at{\Diffeomorphism}{V}$ is
isotopic to the identity as it is the identity on its outer boundary.
Therefore we can assume that it is in fact the identity (note that this changes
the monodromy outside the compression body). Additionally this proves that we
can assume that $\at{\Diffeomorphism}{A}$ is the
identity on the boundary of $A$. Hence
$
  \tau_V
  \coloneqq
  \partial \FundamentalCycle{\Ball{2}}
  \times
  \FundamentalCycle{V}
$
represents the fundamental class of $\Sphere{1} \times V \subset \Manifold$. By
construction the boundary of this cycle agrees with the fundamental cycle
constructed for the tubular neighborhood of the binding on their common
intersection. By choosing $\FundamentalCycle{\Ball{2}}$ sufficiently small,
which results in the norm of $\partial \FundamentalCycle{\Ball{2}}$ being
small, the norm of $\tau_V$ can be made arbitrarily small.

Now consider the components of $\Page \setminus (S\cup V)$, these consist
either of Seifert fibered pieces that can only have boundary components that
are either spheres or tori, or hyperbolic pieces that might have boundary of
genus bigger or equal than $2$. Note that in this case such a boundary
component has to intersect the compression body. Additionally, some of these
boundaries will intersect $A$.

Consider a hyperbolic piece $\Page_{H,i}$ containing a boundary component of
genus greater or equal than $2$.
We will use the following proposition:
\begin{proposition}[Proposition~3.2 in \cite{HatcherMcCullough}]
\label{prp:HatcherMcCullough}
  Let $\Manifold$ denote an irreducible manifold with boundary, such that
  its boundary splits into three pieces intersecting in circles $B \cup A \cup
  T$, where $T$ is the union of all connected components that are tori, $A$ is
  a union of annuli in $\partial \Manifold$ and $B = \partial \Manifold
  \setminus (T\cup A)$ consists solely of components with negative Euler
  characteristic. Assume further that the inclusion of $B$ and $A\cup T$ into
  $\Manifold$ is $\pi_1$-injective, every $\pi_1$-injective map from a torus
  into $\Manifold$ is homotopic into $T$, every map of pairs
  $
    (\Sphere{1}\times [0,1], \Sphere{1}\times \{0,1\})
    \to (\Manifold, B)
  $
  is homotopic as a map of pairs to a map that maps $\Sphere{1}\times [0,1]$
  either into $A$ or $B$ and finally $\Manifold$ is not homeomorphic to an
  interval bundle over a torus.

  Let $R$ denotes a non-empty union of components of $A, B$ and $T$ and let
  $R_0$ denote the components of $R$ with Euler characteristic zero. Then
  $\HomotopyGroupOfObject{0}{\Diff{M,A \text{ rel }R}}{\Identity}$, the mapping
  class group of diffeomorphisms that map $A$ to $A$ and restrict to the
  identity on $R$, is isomorphic to
  $\HomologyOfSpaceObject{1}{R_0}{\Integers}$.
\end{proposition}
Let $R$ denote the union of all components of $\partial \Page_{H,i} \setminus
A$ of negative Euler characteristic. Then applying
\autoref{prp:HatcherMcCullough} yields that $\Diffeomorphism_{H,i}$ is
isotopic to the identity, via an isotopy that is the identity on $R$.
So let us again assume that the monodromy is the identity on these components
(again this changes the monodromy on the page, but not on the compression body!)
Let $\FundamentalCycle{\Page_{H,i}}$ denote a representative
of the fundamental class of $\Page_{H,i}$ which agrees on the corresponding
boundary components with $\FundamentalCycle{V}$ and is also compatible with the
partition of $\partial \Page_{H,i}$ into tori, annuli and the rest.
Then
$
  \tau_{H,i}
  \partial \FundamentalCycle{\Ball{2}} \times \FundamentalCycle{\Page_{H,i}}
$
is a representative of the fundamental class of $\Sphere{1}\times
\Page_{H,i}\subset \Manifold$. Additionally by construction this fundamental
cycle agrees with the so far constructed fundamental cycles on all
intersections of $\Page_{H,i}$ with the inner boundary of $V$ which have
negative Euler characteristic. Again if $\FundamentalCycle{\Ball{2}}$ is chosen
sufficiently small, then the norm of $\tau_{H_i}$ can be made arbitrarily small.

Finally note that by Theorem~VI.18 in \cite{Jaco} any diffeomorphism of a
Seifert fibered piece $\Page_{S,j}$ is isotopic to a fiber-preserving one. Hence
$\MappingTorus{\Diffeomorphism_{S,j}}$ carries an $\Sphere{1}$-action. Hence by
Remark~1 in \cite{Yano} its simplicial volume vanishes. Let
$\tau_{S,j}$ denote a fundamental cycle of
$\MappingTorus{\Diffeomorphism_{S,j}}$ of arbitrarily small norm.

Now
\[
  \tau_0 \coloneqq \tau_\Binding + \tau_V + \sum_i \tau_{H,i} + \sum \tau_{S,j}
\]
almost represents the fundamental class of $\Manifold$, but it is not closed.
Its boundary is not zero.
Nevertheless since the summands of $\tau_0$ represent
the fundamental classes of their corresponding parts of $\Manifold$, their
boundaries represent the fundamental classes of the boundaries of the
corresponding parts.
Now let $\partial_0 V$ denote all spherical components of the inner boundary of
$V$ and let $A'$ denote the union of $A$ with all annular components of
$\partial V \setminus A$.
Then $\partial \tau_0$ can be written as $\rho_V + \rho_S + \rho_A$,
where $\rho_V$ is a difference of two fundamental cycles of
$\MappingTorus{\at{\Diffeomorphism}{\partial_0 V}}$, $\rho_A$ is a difference
of two relative fundamental cycles of
$\MappingTorus{\at{\Diffeomorphism}{\partial_0 A'}}$, which agree on their
boundaries, and $\rho_S$ is a difference of two fundamental cycles of
$\MappingTorus{\at{\Diffeomorphism}{S}}$. Note that the components of
$\MappingTorus{\at{\Diffeomorphism}{\partial_0 V}}$ are mapping tori of
$2$-spheres, the components of $\MappingTorus{\at{\Diffeomorphism}{A'}}$ are
mapping tori of annuli and the components of
$\MappingTorus{\at{\Diffeomorphism}{S}}$ are mapping tori of tori.
Hence they all have an amenable fundamental group.
By Theorem~6.8 in \cite{Frigerio} there exists a $K>0$ such that there exists
chains $\zeta_V$, $\zeta_A$ and $\zeta_S$ supported on
$\MappingTorus{\at{\Diffeomorphism}{\partial_0 V}}$,
$\MappingTorus{\at{\Diffeomorphism}{A'}}$ and
$\MappingTorus{\at{\Diffeomorphism}{S}}$ respectively such that $\partial
\zeta_V = \rho_V$, $\partial
\zeta_A = \rho_A$ and $\partial \zeta_S = \rho_S$ and additionally
\begin{equation}
\label{eqn:UBC}
  \Norm{\zeta_V}
  \leq K
  \Norm{\rho_V}
  \text{and }
  \Norm{\zeta_S}
  \leq K
  \Norm{\rho_S}
\end{equation}
Evidently $\tau_0 + \zeta_V + \zeta_A + \zeta_S$ is closed and represents the
fundamental
class of $\Manifold$. Additionally by (\ref{eqn:UBC}), and the previous
considerations its norm can be chosen to be arbitrarily small.
Hence $\Manifold$ has vanishing simplicial volume, which concludes the proof of
\autoref{prp:Contradiction}.

We have also shown:
\begin{corollary}
  Let $\Manifold$ denote a $3$-manifold and $\Diffeomorphism\colon \Manifold
  \to \Manifold$ a diffeomorphism, which is the identity on the boundary of
  $\Manifold$, then $\MappingTorus{\Diffeomorphism}$ has vanishing simplicial
  volume.
\end{corollary}
Note that it is crucial for this corollary that the diffeomorphism is the
identity on the boundary. Indeed otherwise consider a pseudo Anosov
diffeomorphism in the handlebody group, then the mapping torus of this
diffeomorphism necessarily has non-zero simplicial volume as its boundary is a
hyperbolic $3$-manifold and the simplicial volume of a manifold with boundary
dominates the simplicial volume of the boundary.

\bibliography{sources}

\begin{thebibliography}{McC86}

\bibitem[BN20]{BucherNeofytidis}
Michelle Bucher and Christoforos Neofytidis.
\newblock The simplicial volume of mapping tori of 3-manifolds.
\newblock {\em Math. Ann.}, 376(3-4):1429--1447, 2020.

\bibitem[Bon02]{Bonahon}
Francis Bonahon.
\newblock Geometric structures on 3-manifolds.
\newblock In {\em Handbook of geometric topology}, pages 93--164. Amsterdam:
  Elsevier, 2002.

\bibitem[Fri17]{Frigerio}
Roberto Frigerio.
\newblock {\em Bounded cohomology of discrete groups}, volume 227 of {\em Math.
  Surv. Monogr.}
\newblock Providence, RI: American Mathematical Society (AMS), 2017.

\bibitem[HM97]{HatcherMcCullough}
Allen Hatcher and Darryl McCullough.
\newblock Finiteness of classifying spaces of relative diffeomorphism groups of
  3-manifolds.
\newblock {\em Geom. Topol.}, 1:91--109, 1997.

\bibitem[Jac80]{Jaco}
William Jaco.
\newblock {\em Lectures on three-manifold topology}, volume~43 of {\em Reg.
  Conf. Ser. Math.}
\newblock American Mathematical Society (AMS), Providence, RI, 1980.

\bibitem[McC86]{McCullough}
Darryl McCullough.
\newblock Mappings or reducible 3-manifolds.
\newblock Geometric and algebraic topology, {Banach} {Cent}. {Publ}. 18, 61-76
  (1986)., 1986.

\bibitem[Oer02]{Oertel}
Ulrich Oertel.
\newblock Automorphisms of three-dimensional handlebodies.
\newblock {\em Topology}, 41(2):363--410, 2002.

\bibitem[Qui79]{Quinn}
Frank Quinn.
\newblock Open book decompositions, and the bordism of automorphisms.
\newblock {\em Topology}, 18:55--73, 1979.

\bibitem[Ran10]{Ranicki}
Andrew Ranicki.
\newblock {\em High-dimensional knot theory}.
\newblock Springer monographs in mathematics. Springer, Berlin, Germany,
  December 2010.

\bibitem[Yan82]{Yano}
Koichi Yano.
\newblock Gromov invariant and {{\(S^1\)}}-actions.
\newblock {\em J. Fac. Sci., Univ. Tokyo, Sect. I A}, 29:493--501, 1982.

\end{thebibliography}
\bibliographystyle{alpha}
\end{document}